\documentclass[conference]{IEEEtran}
\IEEEoverridecommandlockouts
\usepackage{cite}
\usepackage{amsmath,amssymb,amsfonts}
\usepackage{algorithmic}
\usepackage{graphicx}
\usepackage{textcomp}
\usepackage{xcolor,mathdots}
\usepackage{hhline}
\usepackage{multirow}
\usepackage{makecell}
\usepackage{amsthm}
\def\BibTeX{{\rm B\kern-.05em{\sc i\kern-.025em b}\kern-.08em
    T\kern-.1667em\lower.7ex\hbox{E}\kern-.125emX}}
\begin{document}

\title{\huge Using covariance extension equation to solve \\the Nevanlinna-Pick interpolation with degree constraint\\
{\footnotesize \textsuperscript{}\large Cui Yufang}}

\maketitle
\begin{abstract}
Nevanlinna-Pick interpolation problem has been widely studied in recent decades, however, the known algorithm is not simplistic and robust enough. This paper provide a new method to solve the Nevanlinna-Pick interpolation problem with degree constraint. It is based on the covariance extension equation proposed by Byrnes and Lindquist. A reformulation of the Nevanlinna-Pick interpolation problem is achieved and then solved by continuation method. This method need not calculate the initial value and a numerical example illustrates robustness and efficiency of the proposed procedure.
\end{abstract}

\begin{IEEEkeywords}
Nevanlinna-Pick interpolation, degree constraint, covariance extension equation, continuation method.
\end{IEEEkeywords}

\section{Introduction}

Nevanlinna-Pick interpolation problem has been widely used in systems and control\cite{b9,b10}, and show obvious advantages in spectral estimation \cite{b2}. To apply this theory better in practice, it is significant to find a robust and efficient algorithm to solve this problem.

In \cite{b8}, Byrnes and Lindquist prove that the space of all strictly positive real, rational functions of at
most degree $ n $ is diffeomorphic to the space of interpolation satisfying the Pick condition and the
space of real Schur polynomials of degree $ n $. 

To solve the problem, a convex optimization method is proposed in \cite{b1}, but the optimization-based procedure in \cite{b1} is numerically infeasible when the interpolant bas poles near the unit circle. Based on the theory in \cite{b8}, a continuation method is applied by constructing a homotopy from the equations for the central solution and the system of interest (see \cite{b3,b4}), but it needs to compute the initial value every time the interpolation points are changed.

This paper describes a new method for the Nevanlinna-Pick interpolation problem with degree constraint. Given $ n+1 $ distinct and self-conjugate points $ z_{0},z_{1},\cdots,z_{n} $ in the complement of the unit disc $ \mathbb{D}^{C}:=\{z||z|>1\} $, and $ n+1 $ self-conjugate values $ w_{0},w_{1},\cdots,w_{n} $ in the open right half-plane. We need to find a rational strict positive real function $ f(z) $ of degree at most $ n $ which satisfies the interpolation condition
\begin{equation}
f(z_{k})=w_{k} \qquad k=0,1,\cdots,n
\end{equation}

The Nevanlinna-Pick theory illustrates that a solution exists if and only if the Pick matrix
\begin{equation}
P_{n}:=\left[ \frac{w_{k}+\bar{w}_{l}}{1-z_{k}^{-1}\bar{z}_{l}^{-1}} \right]_{k,l=0}^{n}
\end{equation}
is positive semidefinite.
This paper is under the following assumption\\
(1) The Pick matrix $ P_{n} $ is positive definite.\\
(2) $ f(z_{i})=w_{i} $ whenever $ f(\bar{z}_{i})=\bar{w}_{i} $\\
(3) $ z_{0}=\infty $, and $ w_{0}$ is real.

The paper is organized as follows. In Section \uppercase\expandafter{\romannumeral2}, we
will review the convex optimization method and show the disadvantages while Section \uppercase\expandafter{\romannumeral3} is devoted to show the basic knowledge of the covariance extension equation(CEE), which is important for the procedure
proposed in this paper. In Section \uppercase\expandafter{\romannumeral4}, we reformulate the Nevanlinna-Pick interpolation problem using CEE. Section \uppercase\expandafter{\romannumeral5} derive a system of nonlinear equation and use a continuation method to solve this equation. An example is given in Section \uppercase\expandafter{\romannumeral6} to illustrate the robustness and efficiency of this method.
\section{The convex optimization solution}

Let
\begin{equation}
\tau (z):=\prod^{n}_{k=1}(z-{z}_{k}^{-1})=z^{n}+\tau_{1}z^{n-1}+\cdots+\tau_{n-1}z+\tau_{n}
\end{equation}
and
\begin{equation}
\rho(z):=\prod^{n}_{k=1}(z-\sigma_{k})=z^{n}+r_{1}z^{n-1}+\cdots+r_{n-1}z+r_{n}
\end{equation}
where $ {z}_{k},k=1,2,\cdots,n $ is the interpolating points and $ \sigma_{k},k=1,2,\cdots,n $ is the spectral zeros.

Let $ \mathcal{H} $ be a subspace consisting all rational, real functions 
\begin{equation}
q(z)=\frac{\pi (z)}{\tau (z)}
\end{equation}
where 
\begin{equation}
\pi (z)=\pi _{0}z^{n}+\pi_{1}z^{n-1}+\cdots+\pi_{n-1}z+\pi_{n}
\end{equation}
for some real numbers $ \pi_{k},k=0,1,\cdots,n $. Then the positive real function $ f(z) $ can be represented by two functions in $ \mathcal{H} $
\begin{equation}
f(z)=\frac{b(z)}{a(z)} 
\end{equation}
where 
\begin{equation}
a(z)=\frac{\alpha (z)}{\tau (z)} \qquad b(z)=\frac{\beta (z)}{\tau (z)}
\end{equation}
Then 
\begin{equation}
\Psi (z):=\frac{\rho(z) \rho(z^{-1})}{\tau(z) \tau(z^{-1})}=a(z)b(z^{-1})+b(z)a(z^{-1})
\end{equation}

Let $ S_{+} $ be the convex set satisfying
\begin{equation}
Q(z):=q(z)+q(z^{-1}) \qquad \text{for some }q\in \mathcal{H}
\end{equation}
and 
\begin{equation}
Q(e^{i\theta})>0 \qquad \theta \in [-\pi,\pi]
\end{equation}

It is clear that $ \Psi(z)\in S_{+} $. For each function $ Q\in S_{+} $, define
\begin{equation}
\begin{split}
J_{\Psi}(Q)=&\frac{1}{2\pi}\int_{-\pi}^{\pi}\{Q(e^{i\theta})[w(e^{i\theta})+w(e^{-i\theta})]\\
&-\log Q(e^{i\theta})\Psi(e^{i\theta})\}d\theta
\end{split}
\end{equation}
where $ w(z) $ is any real function,analytic in $ \{z||z|\geqslant 1\} $ and satisfies $ w(z_{k})=w_{k},\text{for } k=0,1,\cdots,n$

The convex optimization problem 
\begin{equation}
\text{min}\{J_{\Psi}(Q)|Q\in S_{+}\}
\end{equation}
has a unique solution for each $ \Psi\in S_{+} $, see \cite{b1}. Moreover, the unique solution satisfies
\begin{equation}
\frac{\Psi(z)}{Q(z)}=f(z)+f(z^{-1})
\end{equation}
where $ f(z) $ is the positive real function we want to obtain.

To find the unique solution numerically, the Newton's method was proposed in \cite{b1,b2}. It is known in \cite{b1} that the gradient of the function $ J_{\Psi}(Q) $ is infinite at the boundary. Hence, if $ q $ is close to the boundary $ Q_{+} $, this algorithm can be ill-conditioned. Moreover, when the roots of $ \alpha(z) $ near the unit disk, the spectral factorization of $ Q(z) $ is usually hard to solve numerically. 

Some twenty years ago, the covariance extension equation (CEE) was used to solve the rational covariance extension problem(see \cite{b5}), and last year, Lindquist proposed that this nonstandard matrix Riccati equation is universal in the sense that it can be used to solve more general analytic interpolation problems by only changing certain parameters(see \cite{b6}). The following will go deeply into this problem and use CEE to solve the Nevanlinna-Pick interpolation with degree constraint.


\section{The covariance extension equation}

First consider the rational covariance extension problem with degree constraint \cite{b5}. Let
\begin{equation}
c=(c_{0},c_{1},\cdots c_{n})
\end{equation}
a partial covariance sequence such that 
\begin{equation}
T_{n}=\begin{bmatrix}
c_{0}&c_{1}&\cdots&c_{n}\\
c_{1}&c_{0}&\cdots&c_{n-1}\\
\vdots&\vdots&\ddots&\vdots\\
c_{n}&c_{n-1}&\cdots&c_{0}
\end{bmatrix}>0
\end{equation}

For convenience, suppose $ c_{0}=1 $. Given such sequence, the rational covariance extension problem is to find two Schur polynomials
\begin{equation}
a(z)=z^{n}+a_{1}z^{n-1}+\cdots+a_{n}
\end{equation}
\begin{equation}
b(z)=z^{n}+b_{1}z^{n-1}+\cdots+b_{n}
\end{equation}
that satisfy the interpolation condition
\begin{equation}
f(z)=\frac{1}{2}\frac{b(z)}{a(z)}=\frac{1}{2}+c_{1}z^{-1}+\cdots+c_{n}z^{-n}+\cdots
\end{equation}
and the positivity condition
\begin{equation}
f(e^{i\theta})+f(e^{-i\theta})>0,\qquad \text{for all } \theta\in[0,2\pi)
\end{equation}

Then there exists a Schur polynomial
\begin{equation}
\sigma(z)=z^{n}+\sigma_{1}z^{n-1}+\cdots+\sigma_{n}
\end{equation}
satisfying
\begin{equation}
f(z)+f(z^{-1})=w(z)w(z^{-1}):=\Phi(z)
\end{equation}
where
\begin{equation}
w(z)=\rho \frac{\sigma(z)}{a(z)}
\end{equation}
Define
\begin{equation}
\sigma=\begin{bmatrix}
\sigma_{1}\\
\sigma_{2}\\
\vdots\\
\sigma_{n}
\end{bmatrix}\quad 
h=\begin{bmatrix}
1\\
0\\
\vdots\\
0
\end{bmatrix}
\end{equation}
and
\begin{equation}
\Gamma=\begin{bmatrix}
-\sigma_{1}&1&0&\cdots&0\\
-\sigma_{2}&0&1&\cdots&0\\
\vdots&\vdots&\vdots&\ddots&\vdots\\
-\sigma_{n-1}&0&0&\cdots&1\\
-\sigma_{n}&0&0&\cdots&0
\end{bmatrix}
\end{equation}
Furthermore, let $ u_{1},u_{2,\cdots,u_{n}} $ a sequence obtained by expanding 
\begin{equation}
\frac{z^{n}}{z^{n}+c_{1}z^{n-1}+\cdots+c_{n}}=1-u_{1}z^{-1}-u_{2}z^{-1}-\ldots
\end{equation}
and define
\begin{equation}
u=\begin{bmatrix}
u_{1}\\
u_{2}\\
\vdots\\
u_{n}
\end{bmatrix}\quad U=\begin{bmatrix}
0\\
u_{1}&0\\
u_{2}&u_{1}\\
\vdots&\vdots&\ddots\\
u_{n-1}&u_{n-2}&\cdots&u_{1}&0
\end{bmatrix}
\end{equation}
Finally, define the function $ g:\mathbb{R}^{n\times n}\rightarrow \mathbb{R}^{n} $ as
\begin{equation}
g(P)=u+U\sigma+U\Gamma Ph
\end{equation}
Then the Covariance Extension Equation(CEE) has the form
\begin{equation}
P=\Gamma(P-Phh'P)\Gamma'+g(P)g(P)'
\end{equation}
\textbf{Theorem 1:} Let $ (1,c_{1},\cdots,c_{n}) $ be a covariance sequence satisfying (16), Then for each Schur polynomial $ \sigma(z) $, there exists a unique positive semidefinite solution $ P $ of the CEE satisfying $ h'Ph<1 $. Furthermore, there is a unique $ w(z) $ in the form of (23) corresponding to $ \sigma(z) $ and $ (1,c_{1},\cdots,c_{n}) $. The correspondce between $ w(z) $ and $ P $ is
\begin{equation}
\begin{split}
a&=(a_{1},a_{2},\cdots,a_{n})'\\
&=(I-U)(\Gamma Ph+\sigma)-u\\
b&=(b_{1},b_{2},\cdots,b_{n})'\\
&=(I+U)(\Gamma Ph+\sigma)+u\\
\rho&=\sqrt{1-h'Ph}
\end{split}
\end{equation}
This theorem was proved in \cite{b5}.

\section{reformulation of the Nevanlinna-Pick interpolation problem}

Now, consider using CEE to solve the Nevanlinna-Pick interpolation problem \cite{b6}. Since
\begin{equation}
f(z_{k})=\frac{1}{2}\frac{b(z_{k})}{a(z_{k})}=w_{k}
\end{equation}
which can write as 
\begin{equation}
V\begin{bmatrix}
1\\
b
\end{bmatrix}=2WV\begin{bmatrix}
1\\
a
\end{bmatrix}
\end{equation}
where $ V $ is 
\begin{equation}
V=\begin{bmatrix}
z_{0}^{n}&z_{0}^{n-1}&\cdots&1\\
z_{1}^{n}&z_{1}^{n-1}&\cdots&1\\
\vdots&\vdots&~&\vdots\\
z_{n}^{n}&z_{n}^{n-1}&\cdots&1\\
\end{bmatrix}
\end{equation}
and $ W $ is 
\begin{equation}
W=\begin{bmatrix}
w_{0}\\
~&w_{1}\\
~&~&\ddots\\
~&~&~&w_{n}
\end{bmatrix}
\end{equation}
Since the points $ z_{0},z_{1},\cdots,z_{n} $ are distinct, the matrix $ V $ is nonsingular, and hence we have
\begin{equation}
\begin{bmatrix}
1\\
b
\end{bmatrix}=2V^{-1}WV\begin{bmatrix}
1\\
a
\end{bmatrix}
\end{equation} 
Then, by Lemma 5 in \cite{b6},
\begin{equation}
\begin{bmatrix}
0\\
g
\end{bmatrix}=T\begin{bmatrix}
1\\
a
\end{bmatrix}
\end{equation}
where
\begin{equation}
T=\frac{1}{2}(2V^{-1}WV-I)
\end{equation}
finally, according to Lemma 6 in \cite{b6},
\begin{equation}
(I+T)\begin{bmatrix}
0\\
g
\end{bmatrix}=T\begin{bmatrix}
1\\
\Gamma Ph+\sigma
\end{bmatrix}
\end{equation}
Since 
\begin{equation}
I+T=V^{-1}WV+\frac{1}{2}I=V^{-1}(W+\frac{1}{2}I)V
\end{equation}
is nonsingular, define
\begin{equation}
\begin{bmatrix}
u&U
\end{bmatrix}:=\begin{bmatrix}
0&I_{n}
\end{bmatrix}(I_{n+1}+T)^{-1}T
\end{equation}
Then
\begin{equation}
g=u+U\sigma+U\Gamma Ph
\end{equation}
which has the same form as (28), and just change $ u $ and $ U $.
\section{homotopy continuation}

To use the continuation method, we first show the following theorems:\\
\textbf{Theorem 2:}\cite{b8} $ \mathcal{P}_{n} $ is diffeomorphic to $ \mathcal{W}_{n}^{+}\times \delta_{n} $, where $ \mathcal{P}_{n} $ is the space of all pairs of Schur polynomials $ (a,b) $ of the form (17) and (18) such that (19) is strictly positive real, $\mathcal{W}_{n}^{+} $ is the space of all $ w_{1},w_{2},\cdots,w_{n} $ such that (2) is positive definite, and $ \delta_{n} $ is the space of real schur polynomials of degree $ n $. Furthermore, $ (a,b) $ is uniquely determined by $ (w,\sigma)\in(\mathcal{W}_{n}^{+}, \delta_{n}) $.\\
\textbf{Theorem 3:} The map G:$ \mathcal{A}\rightarrow\mathcal{D} $ defined by 
\begin{equation}
G(w,\sigma,P)=(w,\sigma)
\end{equation}
where $ (w,\sigma,P)\in(\mathcal{W}_{n}^{+},\delta_{n},P_{+}) $ that satisfies (29) and (41), is a diffeomorphism. Here, $ P_{+}=\{P\in R^{n\times n}|P\geqslant0,h'Ph<1\} $.

The proof of Theorem 3 is similar to the proof of Theorem 1 in \cite{b7}.\\
\textbf{Proposition 1:} $ w_{k}=1, k=0,\cdots,n $ can make the Pick matrix be positive definite for any given self-conjugate points $ z_{0},z_{1},\cdots,z_{n} $ in the complement of the unit disc which satisfy the assumption $ (2),(3) $\\
\begin{proof}
 It is known that $ f(z) $ can be represented by 
\begin{equation*}
f(z)=\frac{1}{2\pi}\int_{-\pi}^{\pi}\frac{e^{i\theta}+z^{-1}}{e^{i\theta}-z^{-1}}\Phi(e^{i\theta})d\theta,\quad \Phi(e^{i\theta})={\rm Re}\{f(e^{i\theta})\}
\end{equation*}
The problem in this paper can be reformulated as 
\begin{equation}
w_{k}=\frac{1}{2\pi}\int_{-\pi}^{\pi}\frac{e^{i\theta}+z_{k}^{-1}}{e^{i\theta}-z_{k}^{-1}}\Phi(e^{i\theta})d\theta
\end{equation}
Suppose
\begin{equation}
\alpha_{k}(\theta)=\frac{1}{2\pi}\frac{e^{i\theta}+z_{k}^{-1}}{e^{i\theta}-z_{k}^{-1}}
\end{equation}

Since $ z_{0}=\infty $, then $ \alpha_{0} $ is a constant. So, according to \cite{b11}, when we choose $ \Phi=1 $, (43) can generate a sequence such that the Pick matrix is positive definite, precisely,
\begin{equation}
w_{k}=\frac{1}{2\pi}\int_{-\pi}^{\pi}\frac{e^{i\theta}+z_{k}^{-1}}{e^{i\theta}-z_{k}^{-1}}d\theta
\end{equation}
is a positive sequence.

By Residue theorem, we obtain that
\begin{equation}
 w_{k}=1,\quad k=0,1,\cdots,n 
\end{equation}
Then, $ w_{k}=1  $ is a positive sequence 
\end{proof}
To solve the equation (29), we can see that if $ W=\frac{1}{2}I $, which is a positive sequence, then $ u=0,U=0 $, and the equation has the unique solution $ P=0 $. Now consider 
\begin{equation}
\begin{split}
W(\nu)&=\frac{1}{2}I+\nu(W-\frac{1}{2}I),\qquad \nu\in[0,1]\\
&=\begin{bmatrix}
\frac{1}{2}+\nu(w_{0}-\frac{1}{2})\\
~&\ddots\\
~&~&\frac{1}{2}+\nu(w_{n}-\frac{1}{2})
\end{bmatrix}
\end{split}
\end{equation}
and
\begin{equation}
\begin{bmatrix}
u(\nu)&U(\nu)
\end{bmatrix}:=\begin{bmatrix}
0&I_{n}
\end{bmatrix}(I_{n+1}+T(\nu))^{-1}T(\nu)
\end{equation}
where
\begin{equation}
T(\nu)=\frac{1}{2}(2V^{-1}W(\nu)V-I)
\end{equation}
From \cite{b7}, we have
\begin{equation}
ES(a(p))\begin{bmatrix}
1\\
b(p)
\end{bmatrix}=2(1-h'p)d
\end{equation}
where $ p=Ph,E=[I_{n}~0] $
\begin{equation}
d=\begin{bmatrix}
1+\sigma_{1}^{2}+\sigma_{2}^{2}+\cdots+\sigma_{n}^{2}\\
\sigma_{1}+\sigma_{1}\sigma_{2}+\cdots+\sigma_{n-1}\sigma_{n}\\
\sigma_{2}+\sigma_{1}\sigma_{3}+\cdots+\sigma_{n-2}\sigma_{n}\\
\vdots\\
\sigma_{n-1}+\sigma_{1}\sigma_{n}
\end{bmatrix}
\end{equation}
and
\begin{equation}
S(a)=
\begin{bmatrix}

1&\cdots&a_{n-1}&a_{n}\\
a_{1}&\cdots&a_{n}\\
\vdots&\iddots\\
a_{n}
\end{bmatrix}+\begin{bmatrix}
1&a_{1}&\cdots&a_{n}\\
~&1&\cdots&a_{n-1}\\
~&~&\ddots&\vdots\\
~&~&~&1
\end{bmatrix}
\end{equation}
Now, we choose the homotopy $ G: V\times[0,1]\rightarrow R^{n} $ 
\begin{equation}
G(p,\nu):=ES(a(p))\begin{bmatrix}
1\\
b(p)
\end{bmatrix}-2(1-h'p)d
\end{equation}
where $ V:=\{p\in R^{n}|p=Ph,P\in P_{+}\} $

For each $ \nu\in[0,1] $, the system $ G(p,\nu)=0 $ has a unique solution in $ V $, and denote the unique solution as $ p(\nu) $. From Theorem 2 and Theorem 3, the trajectory $ \{p(\nu)\}_{\nu=0}^{1} $ has the following property:

\textbf{Proposition 2:} The trajectory $ \{p(\nu)\}_{\nu=0}^{1} $ is continuously differentiable and has no turning points or
bifurcations.

The above theorems and propositions imply that homotopy continuation method can be used to solve the Nevanlinna-Pick interpolation problem.

The implicit function theorem yields 
\begin{displaymath}
\frac{dp}{d\nu}=\left[\frac{\partial G(p,\nu)}{\partial p}\right]^{-1}\frac{\partial G(p,\nu)}{\partial\nu}
\end{displaymath}
and by calculation, 
\begin{align*}
 \frac{\partial G(p,\nu)}{\partial\nu}  & =ES(a(p,\nu)-b(p,\nu))\begin{bmatrix}0\\\dot{U}(\Gamma p+\sigma) +\dot{u}\end{bmatrix} 
   \\& = -2 ES(U(\Gamma p +\sigma)+u)\begin{bmatrix}0\\\dot{U}(\Gamma p+\sigma) +\dot{u}\end{bmatrix} \\
\frac{\partial G(p,\nu)}{\partial p}  &  =ES(a(p,\nu)-b(p,\nu))\begin{bmatrix}0\\ U\Gamma\end{bmatrix}
\\&+ES(a(p,\nu)+b(p,\nu)))\begin{bmatrix}0\\\Gamma\end{bmatrix} +2hd'
\\&=-2ES(U(\Gamma p +\sigma)+u)\begin{bmatrix}0\\ U\Gamma\end{bmatrix}
\\&+2ES(\Gamma p +\sigma)\begin{bmatrix}0\\\Gamma\end{bmatrix} +2hd'
\end{align*}
where $ u,U,\dot{u},\dot{U} $ are $ u(\nu),U(\nu),\dot{u}(\nu),\dot{U}(\nu) $.
Then we can obtain the following differential equation
\begin{equation}
\begin{split}
&\dot{p}=\left[\frac{\partial G(p,\nu)}{\partial p}\right]^{-1}\frac{\partial G(p,\nu)}{\partial\nu}\\
 &p(0)=0
\end{split}
\end{equation}

This differential equation has a unique solution $ p(\nu) \text{ for every } 0\leq\nu\leq1 $, and the solution of 
\begin{equation}
\begin{split}
&P-\Gamma P\Gamma'=-\Gamma pp'\Gamma'+\\
&\qquad (u+U\sigma+U\Gamma p)(u+U\sigma+U\Gamma p)'
\end{split}
\end{equation}
where $ p=p(1),u=u(1),U=U(1) $, is the unique solution of the equation (29).

To solve above differential equation, we use predictor-corrector steps.
\subsection{Predictor Step}

In the predictor step, Euler's method is used to obtain the estimation $ \hat{p}(\nu+\delta\nu) $ from $ p(\nu) $, more precisely,
\begin{equation}
\hat{p}(\nu+\delta\nu)=p(\nu)+\delta\nu \dot{p}(p(\nu),\nu)
\end{equation}
where $ \dot{p}(p(\nu),\nu) $ is given by (50). The step size $ \delta $ can be determined by the following procedure:
On the trajectory, we have
\begin{equation}
G(p,\nu)=ES(a(p))\begin{bmatrix}
1\\
b(p)
\end{bmatrix}-2(1-h'p)d=0
\end{equation}
and then 
\begin{equation}
e_{1}^{T}G(p,\nu)=e_{1}^{T}(ES(a(p))\begin{bmatrix}
1\\
b(p)
\end{bmatrix}-2(1-h'p)d)=0
\end{equation}
we need that $ \hat{p}(\nu+\delta\nu) $ is close enough to the trajectory, then we can let
\begin{equation}
\mu\leq e_{1}^{T}(ES(a(\hat{p}))\begin{bmatrix}
1\\
b(\hat{p})
\end{bmatrix}-2(1-h'\hat{p})d)\leq\mu
\end{equation}

From the above equation, if $ \mu $ is sufficiently small, $ \hat{p}(\nu+\delta\nu) $ can be close enough to the trajectory.

\subsection{Corrector Step}

In the corrector step, Gauss-Newton algorithm is used to correct the result obtained in the predictor step. The question is to obtain $ p(\nu+\delta\nu) $ such that
\begin{equation}
G(p):=G(p(\nu+\delta\nu),\nu+\delta\nu)=0
\end{equation}

Use $ \hat{p}(\nu+\delta\nu) $ as the initial point and fix $ \nu $ at $ \nu+\delta\nu $, then the iteration formula is
\begin{equation}
\hat{p}_{k+1}=\hat{p}_{k}- \nabla G(\hat{p}_{k})^{-1}G(\hat{p}_{k}),\quad k=0,1,2,\cdots
\end{equation}
where
\begin{align*}
\nabla G(\hat{p}_{k})&=-2ES(U(\Gamma \hat{p}_{k} +\sigma)+u)\begin{bmatrix}0\\ U\Gamma\end{bmatrix}
\\&+2ES(\Gamma \hat{p}_{k} +\sigma)\begin{bmatrix}0\\\Gamma\end{bmatrix} +2hd'
\end{align*}
is invertible when $ \hat{p}_{k} $ is close to the trajectory.
\section{Simulations}

This section will show an example to illustrate the robustness and efficiency compared to the convex optimization method.
\subsection{system identification}
Suppose the interpolation values are
\begin{equation*}
\begin{split}
\{z_{0},\cdots,z_{n}\}=&\{\infty,0.3344-1.2044i,0.3344+1.2044i,\\
&0.8709-0.8967i,0.8709+0.8967i,1.1,\\
&-0.6474-0.8893i,-0.6474+0.8893i\}\\
\{w_{0},\cdots,w_{n}\}=&\{0.5,0.5451 + 0.3645i,0.5451 - 0.3645i,\\
&0.7973+0.2568i,0.7973-0.2568i,0.7693,\\
&0.7693 - 0.7693i,0.7693 + 0.7693i\}
\end{split}
\end{equation*}
which can make the Pick Matrix be positive definite, and suppose the spectral zeros are
\begin{equation*}
\{0.95e^{\pm2.3i},0.95e^{\pm1.22i},\pm0.99i,-0.99\}
\end{equation*}

Since there are 8 interpolation points and 7 spectral zeros in the unit circle, we can obtain a positive real function $ f(z) $ of degree at most 7. After calculation, we find the unique solution $ f(z)=b(z)/2a(z) $, where
\begin{equation*}
\begin{split}
b(z)=&z^{7}-1.364z^{6}+1.112z^{5}-0.3812z^{4}\\
&-0.4479z^{3}+1.119z^{2}-1.412z+0.8781\\
a(z)=&z^{7}-1.771z^{6}+1.815z^{5}-1.205z^{4}\\
&1.28z^{3}-1.814z^{2}+1.773z-0.8775\\
\end{split}
\end{equation*}
and the error
\begin{equation*}
\sum_{t=0}^{5}|w_{j}-f(z_{j})|<10^{-14}
\end{equation*}

At $ \nu=0 $, the zeros of $ a(p(0)) $ are located in Fig.1, and Fig.2 shows the location for the zeros of $ a(p(1))$ or the poles of $ f(z) $  at $ \nu=1 $, Fig.3 shows the trajectories of the poles when vary $ \nu $ from $ 0 $ to $ 1 $.
\begin{figure}[thb!]
  \centering
  \includegraphics[width = 0.8\linewidth]{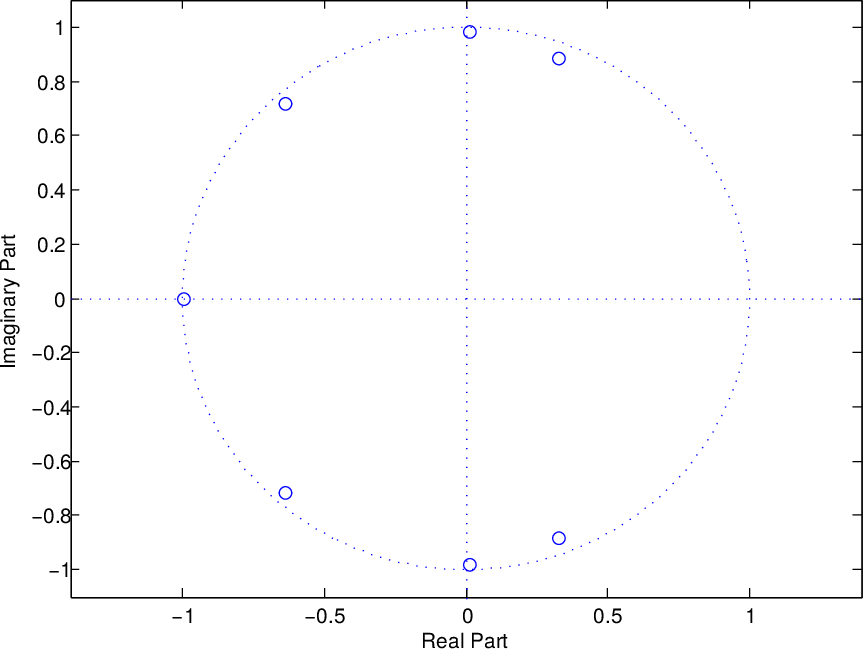}
  \caption{The locations of poles at $ \nu=0 $}
  \label{0}
\end{figure}

\begin{figure}[thb!]
  \centering
  \includegraphics[width = 0.8\linewidth]{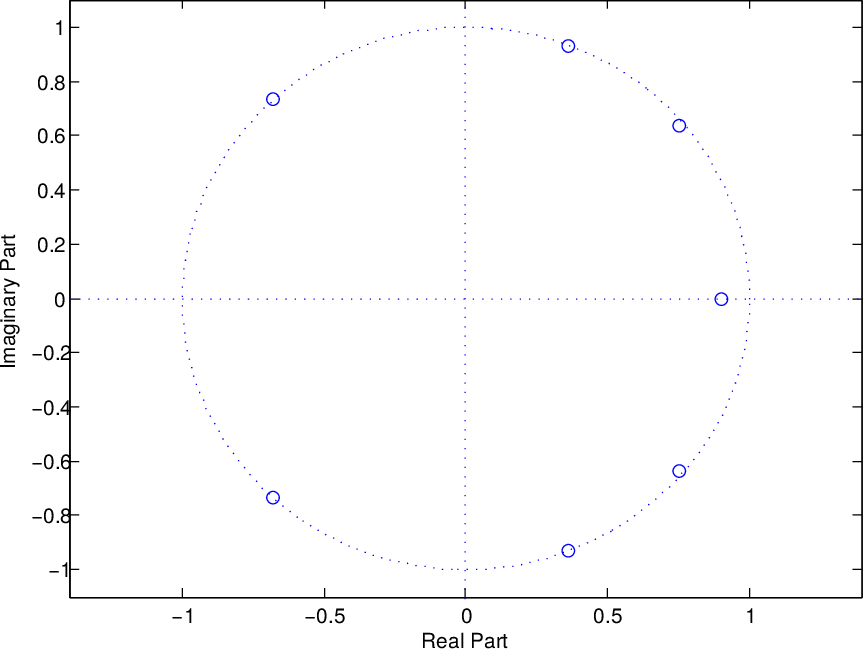}
  \caption{The locations of poles at $ \nu=1 $}
  \label{1}
\end{figure}

\begin{figure}[thb!]
  \centering
  \includegraphics[width = 0.8\linewidth]{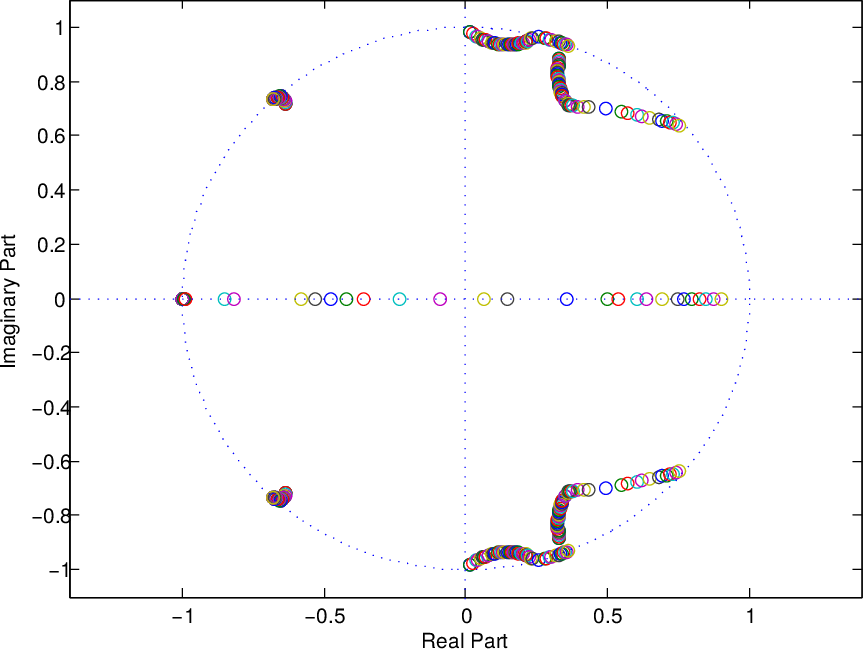}
  \caption{The trajectories of the poles}
  \label{2}
\end{figure}

From the above figures, we can see that some poles of $ f(z) $ are almost on the unit circle, which is hard to solve numerically by the convex optimization method. By contrast, the method in this paper can solve this problem efficiently and robustly. Furthermore, from the trajectories of the poles, the continuity of the map between $\mathcal{W}_{n}^{+} $ and $ \mathcal{P}_{n} $ can easily verified.

in the following, we apply the method described above to detect the positive degree and model reduction. We generate the interpolation data in the following way: Pass white noise through a given stable filter with a rational minimum phase transfer function of degree $ n $,
\begin{equation}
w(z)=\frac{\sigma(z)}{a(z)}
\end{equation}
\begin{figure}[thb!]
  \centering
  \includegraphics[width = 0.7\linewidth]{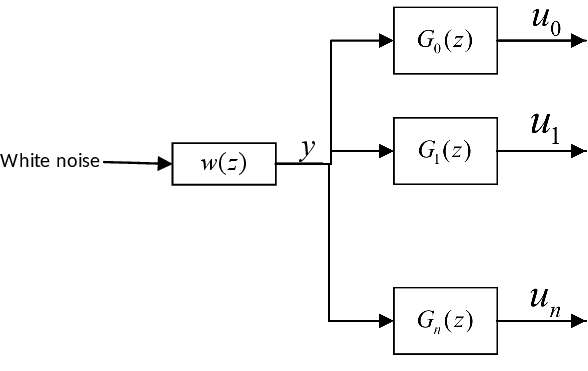}
  \caption{The locations of poles at $ \nu=0 $}
  \label{0}
\end{figure}
then we can obtain the time series $ y_{0},y_{1},y_{2},\cdots,y_{N} $. Pass the output signal to a bank of filters $ G_{0}(z),G_{1}(z),\cdots,G_{n}(z) $, where 
\begin{equation}
G_{k}(z)=\frac{z}{z-p_{k}},|p_{k}|<1\quad k=0,1,\cdots,n
\end{equation}
and the output is $ u_{0},u_{1},\cdots,u_{n} $ respectively. Given any choice of distinct real or complex numbers $ p_{0},p_{1},\cdots,p_{n} $, the values of the positive real function $ f(z) $ at the
points $ p_{0}^{-1},p_{1}^{-1},\cdots,p_{n}^{-1} $ can be expressed by
\begin{equation}
f(p_{k}^{-1})=\frac{1}{2}(1-p_{k}^{2})\mathbb{E}\{u_{k}^{2}\}
\end{equation}
where $ \mathbb{E}\{u_{k}^{2}\} $ can be estimated by
\begin{equation}
\frac{1}{N+1}\sum_{t=0}^{N}u_{k}(t)^{2}
\end{equation}
\subsection{Detecting the positive degree}
Given a transfer function $ w(z) $ of degree $ 2 $ with zeros at $ 0.31e^{\pm 0.98i} $ and
poles at $ 0.76e^{\pm 1.45i} $. The number of interpolation points $ n+1 $ is $ 3,4,5,6,7 $ respectively and set the zero polynomial $ \hat{\sigma}(z)=z^{n-2}\sigma(z) $. Using the method in this paper to compute the $ n\times n $ matrix $ P $. For each $ n $, we perform 100 Monte Carlo simulations. The following table shows the mean value of the singular values of $ P $.

\renewcommand\arraystretch{1.5}
\begin{table}[thb!]
\setlength{\tabcolsep}{0.5mm}{
\caption{Singular values of $ P $}
\begin{tabular}{|c|c|c|c|c|}
\hline
$ n=2 $&$ n=3 $&$ n=4 $&$ n=5 $&$ n=6 $\\ \hline
\makecell{0.2490\\0.2959} & \makecell{0.2494\\0.2968}&\makecell{0.2656\\0.3171}& \makecell{0.2681\\0.3243} & \makecell{0.2731\\0.3448}\\ \hline

\makecell{~\\~\\~\\~} & \makecell{$ 5.0460\times 10^{-6} $\\~\\~\\~}& \makecell{$ 3.1004\times 10^{-6} $\\$ 1.8522\times 10^{-4} $\\~\\~}& \makecell{$ 8.0691\times 10^{-3} $\\$ 3.8549\times 10^{-4} $\\$ 6.3044\times 10^{-5} $\\~}& \makecell{$ 8.5238\times 10^{-3} $\\$ 2.6664\times 10^{-3} $\\$ 2.3724\times 10^{-4} $\\$ 2.5244\times 10^{-5} $}\\ 
\hline
\end{tabular}}
\end{table}

From the table, we can see that for each $ n>2 $, the first two singular values of $ P $ are larger than others which are almost zero, so the rank of $ P $ is close to 2 for $ n\geqslant2 $. Therefore, the positive degree is approximately 2. The following picture shows the given spectral factor for $ n=2 $ and its estimates for $ n>2 $, which almost the same as the given one.

\begin{figure}[thb!]
  \centering
  \includegraphics[width = 0.7\linewidth]{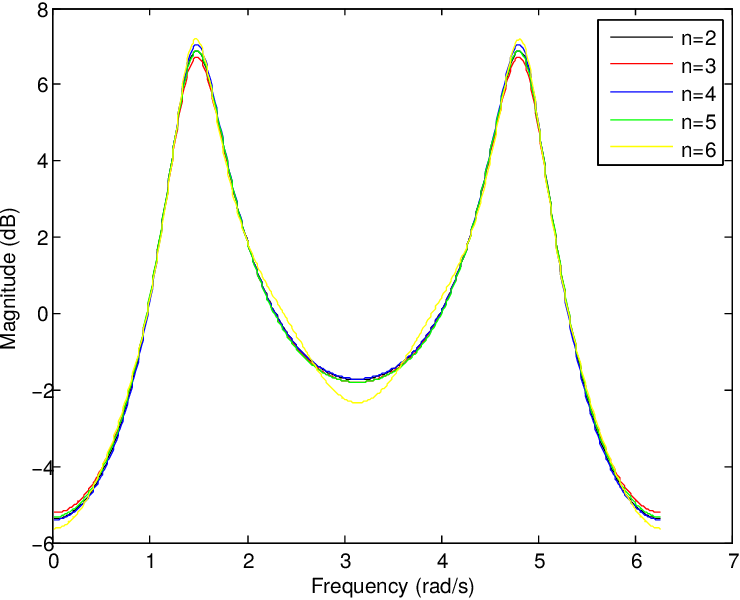}
  \caption{The given spectral factor and its estimated ones}
  \label{0}
\end{figure}

For $ n=4 $, if we change the zero polynomial from $ z^{2}\sigma(z) $ to $ (z-0.6e^{1.5i})(z-0.6e^{-1.5i})\sigma(z) $, then the singular values are $ 0.3399,0.2772,1.3064\times10^{-4},5.5920\times10^{-5} $. Hence the rank of $ P $ is approximately 2. The following figure shows the position of zeros and poles of the estimated spectral factor, which almost has a cancellation at $ 0.6e^{\pm1.5i} $.

\begin{figure}[thb!]
  \centering
  \includegraphics[width = 0.7\linewidth]{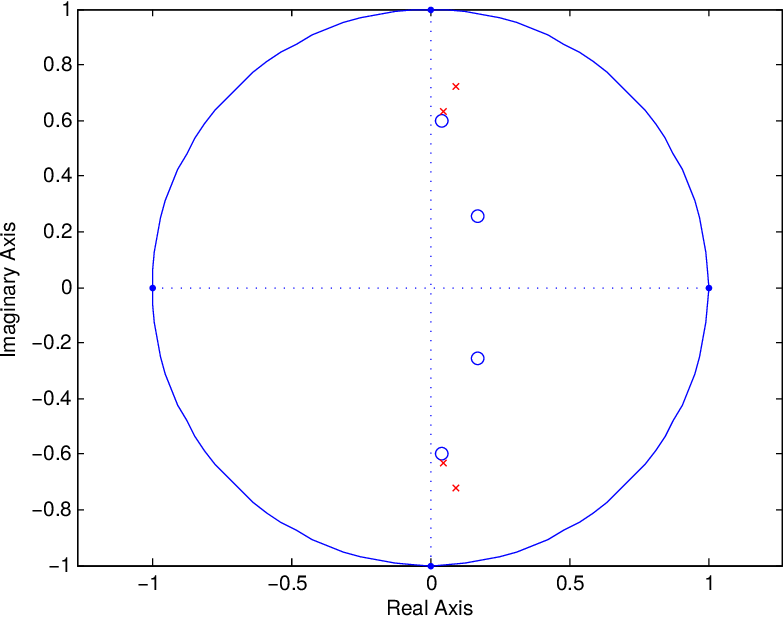}
  \caption{The locations of zeros $ (\circ) $ and poles $ (\times) $ for $ n=4 $}
  \label{0}
\end{figure}

\subsection{Model reduction}

Suppose a transfer function $ w(z) $ of degree 6 with zeros at $ 0.92e^{\pm1.5i},0.49e^{\pm1.4i},0.95e^{\pm2.5i} $ and poles at $ 0.8e^{\pm2.1i},0.83e^{\pm1.34i},0.76e^{\pm0.8i} $. From the following figure, we can see that there is no zero-pole cancellation.
\begin{figure}[thb!]
  \centering
  \includegraphics[width = 0.7\linewidth]{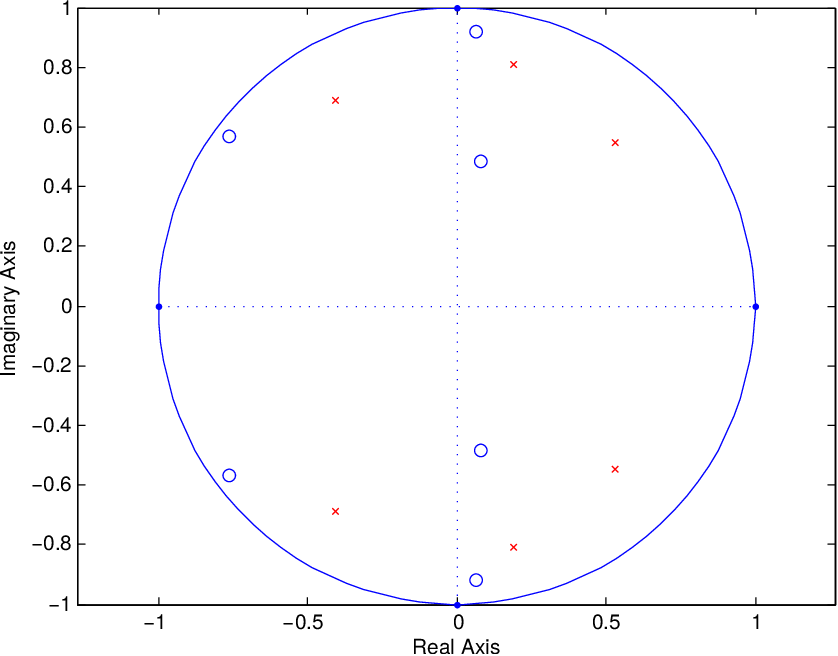}
  \caption{The locations of zeros $ (\circ) $ and poles $ (\times) $ }
  \label{0}
\end{figure}

Using the method in this paper, the singular values of $ P $ is 
\begin{equation*}
2.0170,0.4184,0.02585,0.01858,0.005741,0.002466
\end{equation*}
and the spectral estimation is in figure 8. We can see that the given spectral factor almost the same as the estimated one of degree 6.
\begin{figure}[thb!]
  \centering
  \includegraphics[width = 0.7\linewidth]{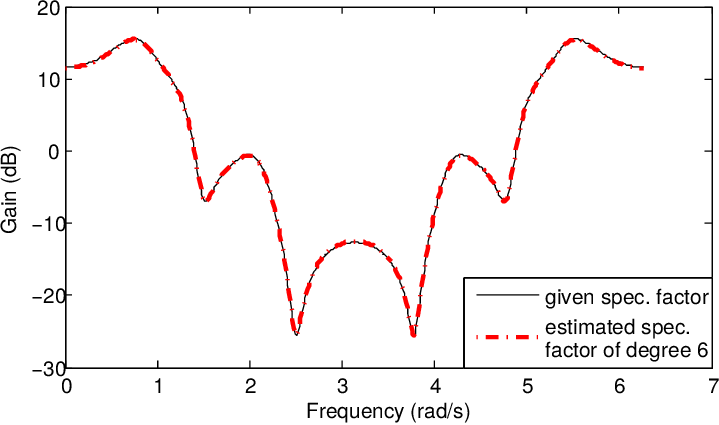}
  \caption{The given spectral factor and its estimated one of degree 6 }
  \label{8}
\end{figure}

Since the last 2 singular values is very small, so the positive degree is close to 4. Now, using the dominant zeros at $ 0.92e^{\pm1.5i},0.95e^{\pm2.5i} $, and then the singular values are
\begin{equation*}
1.2205,0.2913,0.01605,0.02563
\end{equation*} 
and the estimated spectral factor of the reduced order system is in figure 9, showing no much difference. Figure 10 shows the locations of zeros and poles of the reduced order system, and we can see that their locations are quite different.
\begin{figure}[thb!]
  \centering
  \includegraphics[width = 0.7\linewidth]{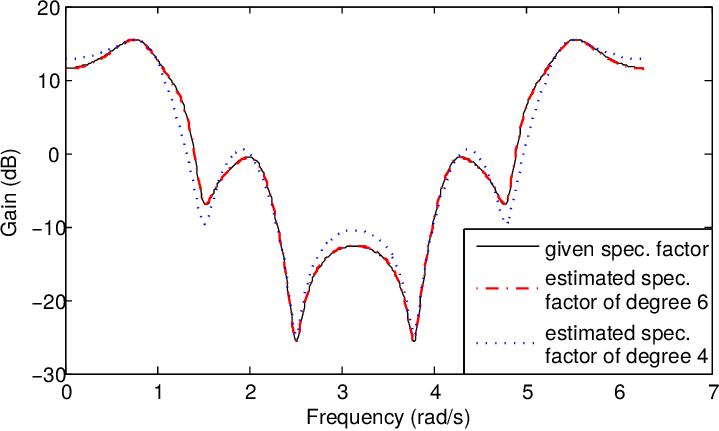}
  \caption{The given spectral factor and its estimated ones}
  \label{9}
\end{figure}

\begin{figure}[thb!]
  \centering
  \includegraphics[width = 0.7\linewidth]{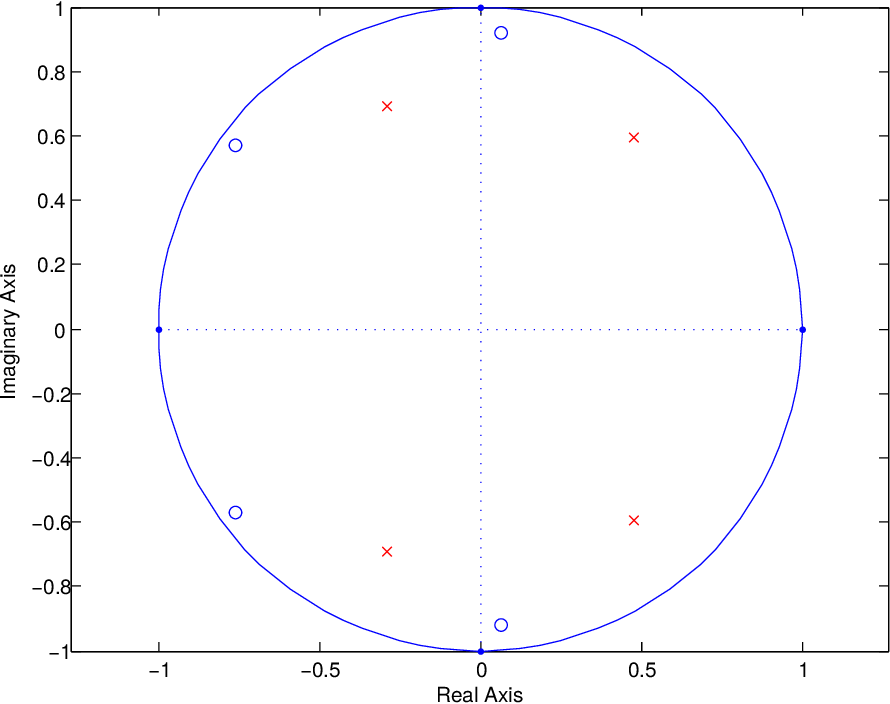}
  \caption{The locations of zeros $ (\circ) $ and poles $ (\times) $ of the reduced order system}
  \label{10}
\end{figure}
\section{conclusion}
This paper provide a new method to solve the Nevanlinna-Pick interpolation problem with degree constraint. The new method reformulates the Nevanlinna-Pick interpolation problem and can be solved by a continuation method with predictor-corrector steps. The numerical experiment shows that the proposed method performs in an robust and efficient manner. This method also verifies that the covariance extension equation can be used to solve Nevanlinna-Pick interpolation problem by only changing certain parameters.

\end{document}